\RequirePackage{cmap} 
\documentclass[12pt]{amsart}
\usepackage[utf8]{inputenc}
\usepackage[T2A]{fontenc}
\usepackage{amsmath,amssymb,amsthm,blkarray,graphicx,mathtools,footmisc, hyperref}
\usepackage[a4paper,hmargin=2.5cm,vmargin=2.5cm]{geometry}
\usepackage[english,russian]{babel}
\usepackage[mathscr]{euscript}
\usepackage{tikz}
\usepackage{tikz-cd}
\usetikzlibrary{arrows.meta,decorations.markings}

\DeclareMathOperator{\coev}{coev}
\DeclareMathOperator{\Rep}{\underline{\smash{\text{Rep}}}}

\tikzset{middlearrow/.style={
        decoration={markings,
            mark= at position 0.55 with {\arrow{#1}} ,
        },
        postaction={decorate}
    }
}

\newcommand{\bbone}{\text{\usefont{U}{bbold}{m}{n}1}}
\MakeRobust{\bbone}

\begin{document}
\author{Ivan Motorin}
\address{Massachusetts Institute of Technology}
\email{ivanm597@mit.edu}
\title[Kac-Moody algebras in Deligne's Category]
      {Kac-Moody algebras in Deligne's Category}

\selectlanguage{english}
\begin{abstract}
  We generalize the notion of a Kac-Moody Lie algebra to the setting of Deligne Categories. Then we derive the Kac-Weyl formula for the category $\mathcal{O}$ integrable representations for such an algebra. This paper generalizes results of A. Pakharev \cite{AP}.
\end{abstract}
\maketitle
\vspace{-0.9cm}

\tableofcontents
\section{Introduction}
The Deligne category $\Rep(GL_t)$ for a parameter $t\in \mathbb{C}$ is a certain interpolation of the representation category of the classical algebraic Lie group $GL_n$ \cite{CW}. One may consider the category $\Rep (\prod_{i\in I} GL_{t_i}):= \boxtimes_{i\in I} \Rep(GL_{t_i})$ of representations of the product group $\prod_{i\in I} GL_{t_i}$. This category is a symmetric tensor category \cite{ERT}. It is possible to define Lie algebra objects or ind-objects within such a category, e.g. a free Lie algebra \cite{EFLie}. Now let $D$ be any graph with the set of vertices \textbf{V}$(D)$ and the set of edges \textbf{E}$(D)$. Given a function $\textbf{N}:\textbf{E}(D)\rightarrow \mathbb{Z}_{\ge 2}$, we may consider the Kac-Moody Lie algebra \cite{Kac} associated to the Dynkin diagram obtained from $D$ by inserting $\textbf{N}(e)-1$ vertices into every edge $e$ from $\textbf{E}(D)$. The goal of this paper is to interpolate the construction of this algebra to the situation when \textbf{N} is replaced by a complex function $\textbf{t}:\textbf{E}(D)\rightarrow \mathbb{C}$ with algebraically independent over $\mathbb{Q}$ values and define it as an ind-algebra within the category $\Rep(\prod_{e\in \textbf{E}(D)} GL_{t_{e}})$ (see Proposition 3.3 and Definition 3.4). The next goal is to study its category $\mathcal{O}$ of representations. This has already been done for classical affine Lie algebras in \cite{AP}.\\

The paper is structured as follows. Chapter 2 contains preliminaries. In Chapter 3 we construct the interpolated Kac-Moody algebra KM$(D,\textbf{t})$. In Chapter 4 we discuss the Lie algebra $\mathfrak{g}$ corresponding to a star-shaped Dynkin diagram consisting of a central vertex with $N$ legs of type $A_{\infty}$, where $N\ge 3$. For the algebra $\mathfrak{g}$ we parameterize $X$-reduced words in the corresponding Weyl group $W$ by particular sequences of nonnegative integers where $X$ is the Dynkin subdiagram of the infinite star-shaped diagram with the central vertex removed (i.e. $N$ disjoint diagrams of type $A_{\infty}$). We also produce an algorithm of descent for these sequences which allows us to algorithmically find the length of the corresponding Weyl words, as we need it for the sign in the Kac-Weyl formula. Finally, in Chapter 5 we define the category $\mathcal{O}$ for KM$(D,\textbf{t})$ by using the structure of $\mathfrak{g}$ and derive the classical Kac-Weyl denominator and character formulas.\\

\textbf{Acknowledgments.} I would like to thank Pavel Etingof for a lot of useful discussions and for supervision of this project. I also thank Aleksei Pakharev for a discussion on his paper. I thank Thomas Lam for his nice argument in Lemma 5.3.
\section{Preliminaries}
\subsection{Kac-Moody Lie algebras.} Suppose we are given an integer square matrix $A$ of size $n$ and rank $l$, such that
\begin{equation}
    a_{ii}=2, \quad a_{ij}\le 0 \text{ if } i\neq j, \quad a_{ij}=0 \Rightarrow a_{ji}=0.
\end{equation}
It is called a generalized Cartan matrix. Let $\mathfrak{h}$ be a vector space of dimension $2n-l$ with independent simple co-roots $\Pi^{\vee}=\{\alpha_1^{\vee},\dots, \alpha_n^{\vee}\}$ in $\mathfrak{h}$ and let $\Pi$ be a set of independent simple roots $\{\alpha_1, \dots, \alpha_n\}$ in $\mathfrak{h}^*$, such that
\begin{equation}
    \langle \alpha_i^{\vee},\alpha_j \rangle =a_{ij}.
\end{equation}
Then there exists a Lie algebra $\mathfrak{g}(A)=\mathfrak{n}_- \oplus \mathfrak{h} \oplus \mathfrak{n}_+$, such that $\mathfrak{n}_+$ is generated by elements $e_1,\dots, e_n$ and $\mathfrak{n}_-$ is generated by elements $f_1,\dots, f_n$ with relations
\begin{equation}
    [e_i,f_j]=\delta_{ij} \alpha_i^{\vee}, \quad [h,h']=0, \quad [h,e_i]=\alpha_i(h)e_i, \quad [h,f_i]=-\alpha_i(h)f_i
\end{equation}
for $h,h'\in \mathfrak{h}$ and
\begin{equation}
    ad_{e_i}^{1-a_{ij}}(e_j)=0, \quad ad_{f_i}^{1-a_{ij}}(f_j)=0.
\end{equation}
Those are called the Chevalley-Serre generators and relations respectively. The constructed Lie algebra is called the Kac-Moody Lie algebra associated to the generalized Cartan matrix $A$ \cite{Kac}.\\

For this Kac-Moody Lie algebra the root lattice $Q$ is the abelian group
\begin{equation}
    Q:=\mathbb{Z}\alpha_1+\dots + \mathbb{Z}\alpha_n.
\end{equation}
The set $Q^+\subset Q$ is the subset of linear combinations of $\alpha_i$ with non-negative coefficients. One can also define the weight lattice $P$ by
\begin{equation}
    P:=\{\lambda\in \mathfrak{h}^*|\forall i \langle \lambda, \alpha_i^{\vee} \rangle \in \mathbb{Z}\}
\end{equation}
and the set of dominant weights $P^+\subset P$,
\begin{equation}
    P^+:=\{\lambda\in \mathfrak{h}^*|\forall i \langle \lambda, \alpha_i^{\vee} \rangle \ge 0\}.
\end{equation}

\subsection{General linear groups.} In this paper we will work over $\mathbb{C}$, but everything here can be generalized to any algebraically closed field $\Bbbk$ of characteristic zero. The general linear group $GL_n(\Bbbk)$ is the group of invertible matrices of size $n$ over $\Bbbk$. The standard choice of a maximal torus $T_n$ of $GL_n(\Bbbk)$ is the subgroup of diagonal matrices and the standard choice of a Borel subgroup $B_n$ is the subgroup of upper-triangular matrices. This yields the following description of the Lie algebra $\mathfrak{gl}_n(\Bbbk)$ of $GL_n(\Bbbk)$ and its root system:
\begin{equation}
    \mathfrak{gl}_n(\Bbbk) = \mathfrak{n}_- \oplus \mathfrak{h} \oplus \mathfrak{n}_+, \quad \mathfrak{n}_-=\text{span}\langle E_{ij} \rangle_{i>j}, \quad  \mathfrak{n}_+=\text{span}\langle E_{ij} \rangle_{i<j}, \quad \mathfrak{h}=\text{span}\langle E_{ii} \rangle,
\end{equation}
\begin{equation}
    \mathfrak{h}^*=\text{span}\langle \epsilon_i \rangle, \quad \epsilon_i(E_{jj})=\delta_{ij}, \quad R=\{\epsilon_{i}-\epsilon_j|i\neq j\}, \quad R^+=\{\epsilon_i -\epsilon_j|i<j\},
\end{equation}
\begin{equation}
    \Pi=\{\epsilon_{i}-\epsilon_{i+1}\},\quad \Pi^{\vee}=\{E_{i,i}-E_{i+1,i+1}\}.
\end{equation}
All irreducible representatios of $GL_n(\Bbbk)$ (or, equivalently, integrable irreducible representations of $\mathfrak{gl}_n(\Bbbk)$) are parameterized by an $n$-tuples of integers $(\lambda_1,\dots,\lambda_n)$ such that $\lambda_i \ge \lambda_{i+1}$. If $V$ is the tautological representation of $GL_n(\Bbbk)$ then any such representation can be tensored with the one-dimensional representation $\Lambda^n V$ several times, so that $\lambda$ becomes a partition of length not greater than $n$. The resulting representation may be realized via a Schur functor $\mathbb{S}^{\lambda}$ applied to $V$.

\subsection{The Deligne category.} The Deligne category $\Rep(GL,T)$ for a formal variable $T$ and a field $\Bbbk$ of characteristic $0$ is the Karoubi closure of the additive closure of the free rigid monoidal $\Bbbk[T]$-linear category generated by an object $V$ of dimension $T$. The endomorphism algebra of an object $V^{\otimes n} \otimes V^{*\otimes m}$ is the walled Brauer algebra $Br_{n,m}(T)$ over $\Bbbk[T]$.\\

For any element $t$ of $\Bbbk$ we may specialize the category $\Rep(GL,T)$ to $T=t$. The resulting $\Bbbk$-linear category $\Rep(GL_t)$ is also usually called a Deligne category \cite{DM}. If $t$ is not an integer then $\Rep(GL_t)$ is abelian and semisimple \cite{CW}. For integer $t$ it is only Karoubian \cite{CW}.\\

Indecomposable objects $L_{[\lambda,\mu]}$ of $\Rep(GL_t)$ are parameterized by bi-partitions $(\lambda,\mu)$ and are obtained by applying appropriate idempotents to $V^{\otimes |\lambda|}\otimes V^{*\otimes |\mu|}$. For any positive integer $t$ the category $\Rep(GL_t)$ admits a full tensor functor $F$ to Rep$(GL_t)$ which sends $V$ to the tautological representation of $GL_t$ and $L_{[\lambda,\mu]}$ to the simple representation in $V^{\otimes |\lambda|}\otimes V^{*\otimes |\mu|}$ with the largest (w.r.t. the standard partial order on the root lattice) highest weight if $l(\lambda)+l(\mu) \le t$. If $l(\lambda)+l(\mu)>t$, then $F(L_{[\lambda,\mu]})=0$.\\

The group $GL_t$ is the fundamental group of $\Rep(GL_t)$ \cite{D}, \cite{ERT}. The Lie algebra $\mathfrak{gl}_t$ (or $\mathfrak{gl}(V)$) of $GL_t$ is
\begin{equation}
    \mathfrak{gl}_t= V\otimes V^*.
\end{equation}
Note that $\mathfrak{gl}_t$ is an associative algebra via the evaluation map, therefore it is also a Lie algebra \cite{ERT}.
\section{Kac-Moody algebras in Deligne's category}
\subsection{Associated Kac-Moody Lie algebras} Consider a graph $D$ with possible multiple edges and self-loops. Let us choose an orientation for each edge $e$, mark the head and tail of the arrow $e$ with the objects $V_e$ and $V_e^*$ in $\Rep(GL_{t_e})$ respectively, and write a parameter $t_e$ on it. For any vertex $v\in \textbf{V}(D)$ we will denote the set of outgoing, incoming and loop arrows for this vertex by $S_v,T_v,L_v$ correspondingly. For any function \textbf{N} as before and $e\in\textbf{E}(D)$ we set $\forall e\in \textbf{E}(D), t_e=\textbf{N}(e)$. Now we consider the same diagram, but with additional $t_e-1$ consecutive inserted vertices which form the type $A_{t_e-1}$ sub-diagram - we will call this the associated diagram.\\
\begin{center}
\begin{tikzpicture}
\filldraw (-3.6,0) circle (2pt);
\filldraw (-1,0) circle (2pt);
\filldraw (1.1,1) circle (2pt);
\filldraw (1.1,-1) circle (2pt);
\draw[thick,->] (-1,0) -- node[below] {$t_1$} (-3.5,0);
\draw[thick,->] (-1,0) -- node[below] {$t_2$} (1,1);
\draw[thick,->] (-1,0) -- node[below] {$t_3$} (1,-1);
\draw[thick,->] (-1,0) to[out=170,in=180] (-1,1) to[out=0,in=45] (-1,0.1);
\node at (-1,1.5) {$t_4$};

\filldraw (4.5,0) circle (2pt);
\filldraw (7,0) circle (2pt);
\filldraw (9,1) circle (2pt);
\filldraw (9,-1) circle (2pt);
\draw[thick] (7,0) -- (4.5,0);
\draw[thick] (7,0) -- (9,1);
\draw[thick] (7,0) -- (9,-1);
\draw[thick] (7,0) to[out=170,in=180] (7,1) to[out=0,in=45] (7,0);
\filldraw (5,0) circle (2pt);
\filldraw (5.5,0) circle (2pt);
\filldraw (6,0) circle (2pt);
\filldraw (6.5,0) circle (2pt);
\filldraw (7+2/3,1/3) circle (2pt);
\filldraw (7+4/3,2/3) circle (2pt);
\filldraw (8,-0.5) circle (2pt);
\filldraw (7,1) circle (2pt);
\end{tikzpicture}\\
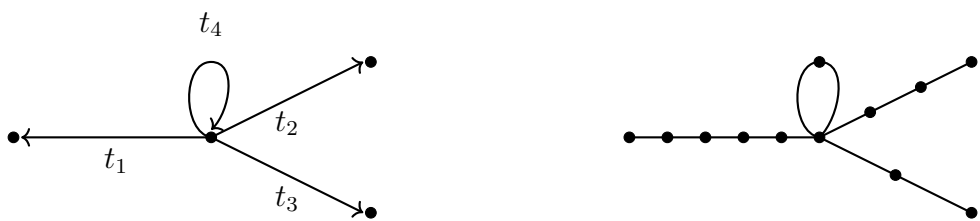
    \textbf{Fig. 1.} Example of a graph $D$ with a self-loop on the left and an associated Dynkin diagram for $t_1=5,t_2=3,t_3=2,t_4=2$.
\end{center}
\vspace{0.5 cm}

Our goal is to interpolate the construction of the simply laced Kac-Moody Lie algebra for the associated diagram, vieweing it as a Lie algebra in the category of representations of the subalgebra $\oplus_{e\in E(D)} \mathfrak{gl}_{t_e}$.\\

Consider the Kac-Moody algebra associated to $D$ and integer $t_e$'s and let us add Cartan elements $Id_e$ corresponding to the edges of the diagram $D$ (the $\mathfrak{gl}_{t_e}$ center). 
Consider an abelian group $\Gamma(D)$ freely generated by vertices $v\in \textbf{V}(D)$. We assign an individual degree to each generator $e_v,f_v$ at a vertex $v\in \textbf{V}(D)$ to be $v$ and $-v$ correspondingly. We will denote such an algebra multi-graded by $\Gamma(D)$ by $\mathfrak{g}(D,\textbf{N})$. We denote a degree $\gamma\in \Gamma(D)$ component of $\mathfrak{g}(D,\textbf{N})$ by $\mathfrak{g}(D,\textbf{N})_{\gamma}$. We also assume that for $Id_e$ in $\mathfrak{gl}_{t_e}$ its action on the Chevalley generator $e_v$ for a vertex $v$ of the edge $e$ is as follows: $[Id_e,e_{v}]=e_v,[Id_e,f_{v}]=-f_v$ for a target vertex $v$ and $[Id_e,e_{v}]=-e_v,[Id_e,f_{v}]=f_v$ for a source vertex $v$. This choice gives us a nice description of $\mathfrak{g}(D,\textbf{N})$.\\

\textbf{Proposition 3.1.} The degree $v$ component of $\mathfrak{g}(D,\textbf{N})$ for a vertex $v\in \textbf{V}(D)$ is
\begin{equation}
    \mathfrak{g}(D,\textbf{N})_v=\bigotimes_{i\in S_v} V_i^* \otimes \bigotimes_{i\in T_v} V_i \otimes \bigotimes_{i\in L_v} L_{i,[(1),(1)]},
\end{equation}
as a $\mathfrak{g}(D,\textbf{N})_0$-module. Here $V_i$ is the tautological $\mathfrak{gl}_{t_i}$-module, $V_i^*$ - its dual and $L_{i,[(1),(1)]}=\ker (ev_{V_i})$:
\begin{equation}
    ev_{V_i}: \mathfrak{gl}_{t_i} = V_i \otimes V_i^* \cong V_i^* \otimes V_i \xrightarrow{ev_{V_i}} \bbone.
\end{equation}
Moreover, the positive part of $\mathfrak{g}(D,\textbf{N})$ is generated by the direct sum of these components of degree $v$. We also have
\begin{equation}
    \mathfrak{g}(D,\textbf{N})_{-v} =\mathfrak{g}(D,\textbf{N})_v^*.
\end{equation}

\textbf{Proof.} The degree $v$ component for $v\in \textbf{V}(D)$ is equal to $U(\mathfrak{g}(D,\textbf{N})_0) e_v$. We may take the positive part of the Lie algebra freely generated by elements $e_v$ for $v\in \textbf{V}(D)$ and for any vertex $y$ adjacent to $v$ impose the Serre relation
\begin{equation}
    ad_{e_{y}}^2(e_v)=0.
\end{equation}
Then the action of the corresponding algebra $U(\mathfrak{gl}_{t_e})$ on $e_v$ will generate one of the modules $V_{e}^*,V_e$ or $L_{e,[(1),(1)]}$ (whether $e\in S_v,T_v$ or $L_v$), because it does so for the associated Dynkin diagram. Here we use the fact that the maximal graded ideal with trivial degree $0$ part in the contragredient Lie algebra without Serre relations for a simply-laced Dynkin diagram is generated by Serre relations (a theorem by O. Gabber, V. Kac  \cite{GK}). Now, since the actions of $\mathfrak{gl}_{t_e}$ commute for different edges $e$, we get that the resulting $U(\mathfrak{g}(D,\textbf{N})_0)$-module is isomorphic to the tensor product of smaller modules as above.\footnote{Note that here we don't have to use the second relation $ad_{e_{v}}^2(e_{y})=0$ due to degree considerations. We also don't use $[e_v,e_w]=0$ for distinct and not connected $v,w\in \textbf{V}(D)$. \label{fnlabel}} The last part follows directly from the above. $\square$\\

\textbf{Remark 1.} We could choose opposite signs for the action of the identity element in $\mathfrak{gl}_{t_e}$. This corresponds to the change of orientation of the edge $e$ and switching of $V_e$ and $V_e^*$.\\

\textbf{Corollary 3.2.} If we don't impose the other Serre relations $ad_{e_{v}}^2(e_{y})=0$ and $[e_v,e_w]=0$ for $v,w\in \textbf{V}(D)$, the positive part of the Lie algebra is isomorphic to the free Lie algebra
\begin{equation} \label{eqn:FLie}
    \text{FLie}(\bigoplus_{v\in \textbf{V}(D)}\bigotimes_{i\in S_v} V_i^* \otimes \bigotimes_{i\in T_v} V_i \otimes \bigotimes_{i\in L_v} L_{i,[(1),(1)]}).
\end{equation}

Now we would like to discuss what happens if we do impose the additional relations $ad_{e_{v}}^2(e_{y})=0$ and $[e_v,e_w]=0$\footref{fnlabel}.\\

\textbf{Proposition 3.3.} The ideal generated by the Serre relations $ad_{e_{v}}^2(e_{y})=0$ and $[e_v,e_w]=0$ in the free Lie algebra \eqref{eqn:FLie} is generated by the direct sum of all of the subspaces below:\\
a)
\begin{equation}
    \bigoplus_{v\in \textbf{V}(D), e\in T_v} L_{e,[(1,1),\o]} \otimes \bigotimes_{i\in S_v} L_{i,[\emptyset, (2)]}\otimes \bigotimes_{i\in T_v,i\neq e} L_{i,[(2), \o]} \otimes \bigotimes_{i\in L_v} L_{i,[(2),(2)]};
\end{equation}
b) 
\begin{equation}
    \bigoplus_{v\in \textbf{V}(D), e\in S_v} L_{e,[\emptyset,(1,1)]} \otimes \bigotimes_{i\in S_v,i\neq e} L_{i,[\emptyset, (2)]}\otimes \bigotimes_{i\in T_v} L_{i,[(2), \o]} \otimes \bigotimes_{i\in L_v} L_{i,[(2),(2)]};
\end{equation}
c) 
\begin{equation}
    \bigoplus_{v\in \textbf{V}(D), e\in L_v} L_{e,[(2),(1,1)]} \otimes \bigotimes_{i\in S_v} L_{i,[\emptyset, (2)]}\otimes \bigotimes_{i\in T_v} L_{i,[(2), \o]} \otimes \bigotimes_{i\in L_v,,i\neq e} L_{i,[(2),(2)]};
\end{equation}
d)
\begin{equation}
    \bigoplus_{v\in \textbf{V}(D), e\in L_v} L_{e,[(1,1),(2)]} \otimes \bigotimes_{i\in S_v} L_{i,[\emptyset, (2)]}\otimes \bigotimes_{i\in T_v} L_{i,[(2), \o]} \otimes \bigotimes_{i\in L_v,,i\neq e} L_{i,[(2),(2)]};
\end{equation}
e)
\begin{equation}
    \bigoplus_{e\in \textbf{E}(D)\setminus  \textbf{L}(D)} L_{e,[(1),(1)]} \otimes \bigotimes_{i\in S_{s(e)} \cup S_{t(e)} ,i\neq s(e)} V_i^* \otimes \bigotimes_{i\in T_{s(e)} \cup T_{t(e)},i\neq t(e)} V_i \otimes \bigotimes_{i\in L_{s(e)} \cup L_{t(e)}} L_{i,[(1),(1)]};
\end{equation}
for different $v,w\in \textbf{V}(D)$ which are not connected by an edge\\
f) 
\begin{equation}
    \bigotimes_{i\in S_v \cup S_w} V_i^* \otimes \bigotimes_{i\in T_v\cup T_w} V_i \otimes \bigotimes_{i\in L_v\cup L_w} L_{i,[(1),(1)]}.
\end{equation}
Here $s(e)$ and $t(e)$ denote the source and the target of an edge $e$ correspondingly, and $\textbf{L}(D)$ is the set of loop edges for $D$.\\

\textbf{Proof.} Recall that the commutator of a free Lie algebra FLie$(V)$ in degree $2$ is the anti-symmetrization map
\begin{equation}
    [,]: V \otimes V \rightarrow \Lambda^2 V.
\end{equation}
Also note that for vector spaces $V_i,i\in I$ labeled by an ordered set $I$
\begin{equation}
    \Lambda^2 (\oplus_{i\in I} V_i) = \bigoplus_{i\in I} \Lambda^2 V_i \oplus \bigoplus_{i,j\in I,i< j} V_i \otimes V_j
\end{equation}
and for finite $I$
\begin{equation}\label{eq:TL2}
    \Lambda^2(\otimes_{i\in I}V_{i}) = \bigoplus_{S\subset I,|S| \text{ is odd}} \bigotimes_{i\in S} \Lambda^2 V_i \otimes \bigotimes_{i\in I\setminus S} S^2 V_i.
\end{equation}
Similarly
\begin{equation}
    S^2(\oplus_{i\in I} V_i) =\bigoplus_{i\in I} S^2 V_i \oplus \bigoplus_{i,j\in I,i< j} V_i \otimes V_j
\end{equation}
and
\begin{equation}\label{eq:TS2}
    S^2(\otimes_{i\in I}V_{i}) = \bigoplus_{S\subset I,|S| \text{ is even}} \bigotimes_{i\in S} \Lambda^2 V_i \otimes \bigotimes_{i\in I\setminus S} S^2 V_i.
\end{equation}
Therefore, we may decompose the second degree part of the free Lie algebra \eqref{eqn:FLie} into the simple $\mathfrak{g}(D,\textbf{N})_0$-modules and analyze where the Serre relations of the form $[e_v,[e_v,e_i]]$ or $[e_v,e_w]$ lie (here $e_i$ is the closest to $v$ generator from an edge adjacent to $v$). It is easy to see that the adjoint action of $e_i$ on $e_v$ translates to the action only on the highest vector of the corresponding simple module $V_i$ (here we abuse the notation by referring to the adjacent to $v$ vertex $i$ on the corresponding edge). That's why we have only one exterior power and the highest weight simple submodules in the symmetric powers for the other modules in the cases a) and b).\\

In the cases c) and d) we have two similar relations but for the same edge. It's not hard to see that for $\mathfrak{gl}_n$ with standard basis $E_{ij}$, a choice of a Borel subalgebra and dual Cartan basis $\epsilon_i$ we have that
\begin{equation}
    E_{n,1} \wedge [E_{1,2},E_{n,1}] \text{ and } E_{n,1} \wedge [E_{n-1,n},E_{n,1}]
\end{equation}
are both lowest weight vectors of weights $-\epsilon_1-\epsilon_2+2\epsilon_n$ and $-2\epsilon_1+ \epsilon_{n-1}+\epsilon_n$ correspondingly. Then the highest weight vectors of such submodules will have weights $\epsilon_1+\epsilon_2-2\epsilon_n$ and $2\epsilon_1- \epsilon_{n-1}-\epsilon_n$ correspondingly.\\

It is useful to provide a decomposition of $\Lambda^2(V\otimes V^*)$ and $S^2(V\otimes V^*)$ into simple $V\otimes V^*$-modules. According to \eqref{eq:TL2} and \eqref{eq:TS2} we have
\begin{equation}\label{eq:crutch1}
    \Lambda^2(V\otimes V^*) \cong \Lambda^2 V \otimes S^2 V^* \oplus S^2 V \otimes \Lambda^2 V^*,
\end{equation}
\begin{equation}\label{eq:crutch2}
    S^2(V\otimes V^*) \cong \Lambda^2 V \otimes \Lambda^2 V^* \oplus S^2 V \otimes S^2 V^*.
\end{equation}
Furthermore, the contraction $ev_{1,3}$ along the $1$st and $3$rd tensor components gives us a $V\otimes V^*$-invariant map from both of the summands in \eqref{eq:crutch1} and \eqref{eq:crutch2} into $\ker(ev_v)$ (i.e. $\mathfrak{sl}(V)$), so we have
\begin{equation}
    \Lambda^2(V\otimes V^*) \cong L_{[(1,1),(2)]}\oplus L_{[(2),(1,1)]} \oplus L_{[(1),(1)]]}\oplus L_{[(1),(1)]]}\oplus ?,
\end{equation}
but the dimension comparison (using the Weyl dimension formula) tells us that $?$ is zero.\\

Analogously, 
\begin{equation}
    S^2(V\otimes V^*) \cong L_{[(2),(2)]}\oplus L_{[(1,1),(1,1)]} \oplus L_{[1,1]}^{\oplus 2} \oplus L_{\emptyset,\o}^{\oplus 2},
\end{equation}
because the contraction map $ev_{13}$ now maps both $S^2V\otimes S^2 V^*$ and $\Lambda^2V\otimes \Lambda^2 V^*$ surjectively onto $V\otimes V^* \cong L_{[1,1]} \oplus L_{\emptyset,\o}$. The rest follows from the Weyl dimension formula.\\

Finally, in the cases e) and f) of two different vertices $v,w$ the quadratic relation lies in the tensor product of the corresponding degree one components and it's not hard to see that it generates the suggested submodules. $\square$\\

\subsection{Categorical Kac-Moody algebra.} The discussion above motivates the following definition.\\

\textbf{Definition 3.4.} Consider the category $\Rep(\prod_{e\in E(D)}GL_{t_e})$ for algebraically independent transcendental parameters $t_e\in \mathbb{C}$. Define the Kac-Moody algebra KM$(D,\textbf{t}),\textbf{t}(e)=t_e$ in this category as an ind-object equal to the direct sum of its Levi subalgebra $\mathfrak{l}$ and nilpotent subalgebras $\mathfrak{u}_+,\mathfrak{u}_-$
\begin{equation}
    \text{KM}(D,\textbf{t}) = \mathfrak{u}_- \oplus \mathfrak{l} \oplus \mathfrak{u}_+
\end{equation}
\begin{equation}
    \mathfrak{l}= \bigoplus_{e\in \textbf{E}(D)} V_{e}\otimes V_{e}^* \oplus \bigoplus_{v\in \textbf{V}(D)} (\bbone_v\oplus \bbone_v'). 
\end{equation}

$\bullet$ The positive part $\mathfrak{u}_+$ is equal to the quotient of \eqref{eqn:FLie} by the ideal generated by the objects in Proposition 3.3 and the negative part $\mathfrak{u}_-$ is equal to the dual free Lie algebra modulo the dual ideal.\\

$\bullet$ The adjoint action of $\mathfrak{gl}(V_e)\subset \mathfrak{l}$ for an edge $e\in\textbf{E}(D)$ in $\mathfrak{l}$ on $\mathfrak{u}_+$ or $\mathfrak{u}_-$ is obvious. For $v\in \textbf{V}(D)$ the action of $\bbone_{v}$ corresponds to the vertex multi-grading by $\Gamma(D)$
and the action of $\bbone_v'$ is zero. \\

$\bullet$ The commutator between degree $v$ and $-w$ parts of $\mathfrak{u}_+$ and $\mathfrak{u}_-$ for distinct $v,w\in \textbf{V}(D)$ is zero. And the commutator between degree $v$ and $-v$ parts of $\mathfrak{u}_+$ and $\mathfrak{u}_-$ correspondingly is given by the formula 
\begin{multline}
    [,]: \bigotimes_{i\in S_v} V_i^* \otimes \bigotimes_{i\in T_v} V_i \otimes \bigotimes_{i\in L_v} L_{i,[(1),(1)]} \otimes \bigotimes_{i\in S_v} V_i \otimes \bigotimes_{i\in T_v} V_i^* \otimes \bigotimes_{i\in L_v} L_{i,[(1),(1)]} \rightarrow\\
    \rightarrow \bigoplus_{i\in S_v \cup T_v \cup L_v} V_i \otimes V_i^* \oplus  \bbone_v\oplus \bbone_v', 
\end{multline}
where the map to each component $V_i\otimes V_i^*$ is a composition of contractions along irrelevant factors (the contraction between two $L_{l,[(1),(1)]}$'s is the trace form) and the maps
\begin{equation}
    V_i \otimes V_i^* \xrightarrow{Id} V_i \otimes V_i^*,\quad V_i^* \otimes V_i \xrightarrow{-(12)} V_i \otimes V_i^*,    
\end{equation}
\begin{equation}
    \quad L_{l,[(1),(1)]} \otimes L_{l,[(1),(1)]} =\mathfrak{sl}(V_l)\otimes \mathfrak{sl}(V_l)\xrightarrow{[,]_{\mathfrak{gl}_{t_l}}} V_l \otimes V_l^*.
\end{equation}
The maps to $\bbone_v$ and to $\bbone_v'$ are the full contraction times $2-|T_v|-|S_v|-2|L_v|$ and the full contraction correspondingly. \\ 

$\bullet$ The full commutator map in $\text{KM}(D,\textbf{t})$ is determined from this data by the Jacobi identity.\\

\textbf{Remark 2.} This construction does not reflect the situation when one of the parameters $t_e$ is equal to $1$, because in that case we would need to impose two cubic relations between generators of adjacent vertices instead of a quadratic one.\\

\textbf{Proposition 3.5.} KM$(D,\textbf{t})$ is a Lie algebra for algebraically independent transcendental $t_e$.\\

\textbf{Proof.} The only thing we have to show is that for any three irreducible subobjects $X,Y,Z$ of KM$(D,\textbf{t})$ the map
\begin{equation}
    [*,[*,*]]\circ (id+ (123)+(132)): X\otimes Y\otimes Z \rightarrow \text{KM}(D,\textbf{t})
\end{equation}
is identically zero. For integer values of our parameters $t_e$ we have a standard symmetric tensor functor $F$ from $\Rep(\prod_{e\in \textbf{E}(D)} GL_{t_{e}})$ to Rep$(\prod_{e\in \textbf{E}(D)} GL_{t_{e}})$ which for a fixed graded part of KM$(D,\textbf{t})$ and large enough parameters $t_e$ sends it to the respective graded part of the associated Kac-Moody Lie algebra with some additional center. This functor is fully faithful and surjective there in the sense that the decomposition of these components in both algebras involves the same number (and type) of irreducible $\mathfrak{l}$-modules. Clearly, then
\begin{equation}\label{eq:Jacobi}
    F([*,[*,*]]\circ (id+ (123)+(132)))=0.
\end{equation}

Now, suppose that $X,Y,Z$ and images of each term of the map above are simultaneously well-defined objects in the classical representation category $\Rep(\prod_{e\in E(D)} GL_{t_{e}})$ (i.e. every simple term appearing in the $K$-group decomposition is non-zero) for some integer parameters $t_e$. Then it's not hard to see that the map \eqref{eq:Jacobi} can be realized as a morphism in a bigger space expressible in a uniform way via $(w,w')$-diagrams with coefficients in $\mathbb{Q}(t_e)$ \cite{CW}. But if the map above becomes zero under $F$-specialization for any large enough integers $t_e$, then it had to be zero to begin with. $\square$\\

\textbf{Remark 3.} This proof reflects the essence of the interpolation argument for Deligne categories. The result is also true for not algebraically independent and transcendental $t_e$, but then one has to use interpolation of the Deligne category $\Rep(GL_t)$ from a case of positive characteristic (see Theorem 3.11 and Lemma 3.12 in \cite{EK}). For integer $t_e$ this construction allows us to span the whole co-root Cartan subspace exactly as in the Kac-Moody Lie algebra. The addition of objects $\bbone_v'$ in necessary only for vertices of valency $2$, as it guarantees the linear independence of co-roots in the Cartan subalgebra.\\

\textbf{Remark 4.} For variables $t_e$ which are not algebraically independent and transcendental we can define the Kac-Moody Lie algebra in three ways. Namely, let $I_1$ be the graded ideal generated by Serre relations, $I_2$ be a limit (which might depend on the direction from which we approach $t_e$) graded ideal interpolated by ideals generated by subobjects from Proposition 3.3 and $I_3$ be the radical of the invariant form. We have the chain of embeddings $I_1 \subset I_2 \subset I_3$. Then the Kac-Moody algebra should be defined as a quotient of the corresponding contragredient Lie algebra by $I_3$. We have proved that for generic $t_e$ the ideals $I_1$ and $I_3$ coincide, but it is unknown to us whether $I_1$ is equal to $I_3$ for special $t_e$.\\

\textbf{Remark 5.} If the graph $D$ has a leaf vertex $v$ with the adjacent edge $e$, then one may drop all nontrivially $v$-graded subobjects in $\text{KM}(D,\textbf{t})$  and $\bbone_v,\bbone_v'$. The resulting algebra $(\text{KM}(D,\textbf{t}),v)$ corresponds to an amputated graph $D$ with omitted vertex $v$. Furthermore, if we apply the restriction functor $\text{Res}_e: \Rep(GL_{t_e}) \rightarrow \Rep(GL_{t_e-1})$ (extended to the whole product category) to $(\text{KM}(D,\textbf{t}),v)$, we will obtain the original algebra $\text{KM}(D,\textbf{t}')$ without $\bbone_v,\bbone_v'$ for the same function $\textbf{t}'$ as $\textbf{t}$ except for $\textbf{t}'(e)=\textbf{t}(e)-1$. This can also be done for a set of leaf vertices.\\

Now we would like to describe the center of KM$(D,\textbf{t})$.\\

\textbf{Proposition 3.6.} The images of the following maps generate the center of KM$(D,\textbf{t})$:
\begin{equation}
    (1,\coev_{V_e},-1): \bbone \rightarrow \bbone_{s(e)}\oplus V_e \otimes V_e^* \oplus \bbone_{t(e)},\quad e\in \textbf{E}(D)\setminus \textbf{L}(D);
\end{equation}
\begin{equation}
    \coev_{V_{l}}: \bbone \rightarrow  V_l \otimes V_l^*, l\in \textbf{L}(D);
\end{equation}
and
\begin{equation}
    1_v': \bbone \rightarrow \bbone_v', v\in \textbf{V}(D).
\end{equation}

\textbf{Proof:} Obviously, the images above lie in the center. For the converse, note that any other central map from $\bbone$ to $\mathfrak{l}\subset \text{KM}(D,\textbf{t})$ (remember about grading by $\Gamma(D)$) should not have any nontrivial components in $L_{e,[(1),(1)]} \subset \mathfrak{l}$ for any $e\in \textbf{E}(D)$. Without loss of generality we may assume that the central map has nontrivial components only at $\bbone_v, v\in \textbf{V}(D)$, but due to the nontrivial action of these objects the central map must be identically zero. $\square$
\section{Star-shaped diagrams}

In this section we will study the Lie algebra with a star-shaped Dynkin diagram with $N$ legs of type $A_{\infty}$. Consider an amputated graph $D$ with one vertex $v$ and $N$ arrows oriented towards $v$, then the corresponding Kac-Moody Lie algebra can be described as
\begin{equation}
    \mathfrak{g} = (\text{FLie}(\bigotimes_{i=1}^N V_i^*)/I^*)\oplus \bigoplus_{i=1}^N \mathfrak{gl}(V_i) \oplus \bbone_v \oplus \bbone_v' \oplus(\text{FLie}(\bigotimes_{i=1}^N V_i)/I)
\end{equation}
where $I$ is the ideal in the free Lie algebra generated by
\begin{equation}
    \bigoplus_{i=1}^N (L_{i,[(1,1),\o]} \otimes \bigotimes_{j\neq i}L_{j,[(2),(0)]})
\end{equation}
and $I^*$ is the ideal generated by dual relations. Sometimes we will omit $\bbone_v'$.\\

This algebra can be viewed as a colimit of corresponding algebras with finite legs. The case of $N=1$ has an obvious description via an amputated segment graph and the case $N=2$ was already considered in \cite{ERT} and in \cite{AP}, because this is the case of $\mathfrak{gl}_{t+s}$ restricted to $\mathfrak{gl}_{t}\oplus \mathfrak{gl}_s$. So assume that $N\ge 3$ and orient the edges towards the central vertex.\\

\subsection{The Lie algebra} Let us fix some notation. The Lie algebra $\mathfrak{g}$ associated with the aforementioned diagram has decomposition $\mathfrak{n}_-\oplus \mathfrak{h} \oplus \mathfrak{n}_+$ with $\mathfrak{n}_+$ generated by elements $E_{i+1,i}^k\in \mathfrak{gl}_{\infty,k},k\in \overline{1,N},i\in \mathbb{N}$ and $E$ and the negative part is generated by $E_{i,i+1}^k\in \mathfrak{gl}_{\infty,k},k\in \overline{1,N},i\in \mathbb{N}$ and $F$. The Cartan subalgebra $\mathfrak{h}$ is spanned by $E_{ii}^k,k\in \overline{1,N},i\in \mathbb{N}$ and $1_v$. This algebra also has the Levi decomposition
\begin{equation}
    \mathfrak{g}=\mathfrak{u}_-\oplus \mathfrak{l}  \oplus \mathfrak{u}_+,\quad \mathfrak{l}=\bigoplus_{k=1}^N \mathfrak{gl}_{\infty,k}\oplus \bbone_v 
\end{equation}
with commutator relations as in the previous chapter. Note that we do not include the $\bbone_v'$ space here, because in the case $N\neq 2$ it is purely superfluous.\\

Consider the dual space $\mathfrak{h}^*$ of finitely supported (w.r.t. the chosen basis of $\mathfrak{h}$) linear functions on $\mathfrak{h}$ with dual basis $\epsilon_{i,k},\epsilon_{v}$.
\begin{equation}
    \epsilon_{i,k}(E_{n,n}^m)=\delta_{i,n}\delta_{k,m},\quad \epsilon_{i,k}(1_v)=\epsilon_v(E_{n,n}^m)=0, \quad \epsilon_v(1_v)=1
\end{equation}
In the extended case we can also add orthogonal symbols
\begin{equation}
    \epsilon_v'\in \mathfrak{h}^*,\quad 1_v'\in \mathfrak{h},\quad \epsilon_v'(1_v')=0.
\end{equation}
Then the simple positive roots will be of the following form:
\begin{equation}
    \alpha_{i,k}=\epsilon_{i+1,k}-\epsilon_{i,k},\quad \alpha_v= \epsilon_v+\sum_{k=1}^N\epsilon_{1,k}
\end{equation}
and the co-roots will be of the form 
\begin{equation}
    E_{i+1,i+1}^k-E_{i,i}^k,\quad (2-N)1_v+\sum_{k=1}^N E_{1,1}^k
\end{equation}
or
\begin{equation}
    1_v'+(2-N)1_v+\sum_{k=1}^N E_{1,1}^k
\end{equation}
if we do want to include the action of $\bbone_v'$. The reflection $s_{i,k}$ w.r.t. the root $\epsilon_{i+1,k}-\epsilon_{i,k}$ transposes $\epsilon_{i,k}$ and $\epsilon_{i+1,k}$ and doesn't change $\epsilon_v$. The reflection $s_v$ with respect to the central root acts nontrivially only on $\epsilon_v,\epsilon_{1,k}$:
\begin{equation}
    s_v(\epsilon_v)=(N-1)\epsilon_v +(N-2)\sum_{k=1}^N \epsilon_{1,k}, \quad s_v(\epsilon_{1,k})=-\epsilon_v -\sum_{j\neq k} \epsilon_{1,j}
\end{equation}
It also acts nontrivially on $\epsilon_v'$:
\begin{equation}
    s_v(\epsilon_v')=\epsilon_v' -\epsilon_v - \sum_{k=1}^N \epsilon_{1,k}.
\end{equation}
Let $a_{v},a_{ik}$ be the dual basis to $\epsilon_v,\epsilon_{1,k}$ in $(\mathfrak{h}^*)^*$, then the generators of the Weyl group act by transpositions for $s_{i,k}$ and
\begin{equation}
    s_v(a_v)=(N-1)a_v -\sum_{k=1}^N a_{1,k},\quad s_v(a_{1,k})=(N-2)a_v -\sum_{j\neq k} a_{1,j}.
\end{equation}

In order to motivate the next Proposition, we prove the following.\\

\textbf{Proposition 4.1.} Consider the ring of formal power series $\mathbb{C}[[a_v,a_{i,k}]]$. Then the invariants of this ring under the action of the Weyl group $W$ generated by $s_v,s_{i,k}$ are generated by
\begin{equation}\label{eq:Inv}
    -(N-2)a_v^2 +\sum_{i,k}a_{i,k}^2 \text{ and } a_v-\sum_{i} a_{i,k}, \forall k.
\end{equation}

\textbf{Proof.} We may identify $\mathfrak{h}$ and $\mathfrak{h}^*$ via the quadratic function in \eqref{eq:Inv}. One can check that the functions above are indeed $W$-invariant. For the converse, we may assume that a $W$-invariant function is homogeneous with respect to $a_v,a_{i,k}$. Let us pick $M\in \mathbb{N}$ and let $a_{i,k}=0$ for $i>M$. Note then that the computation of invariant polynomials in this case is equivalent to the computation of invariants of the Cartan polynomial algebra $\mathfrak{h}^{W_M}_M$ of a direct sum of the Kac-Moody Lie algebra corresponding to the star-shaped diagram with $N$ legs of length $M-1$ with $N$-dimensional center spanned by $1_v-\sum_{i=1}^{M}E_{i,i}^k$. Since the $W$-action on $a_v,a_{i,k}$ has real coefficients, we know that the invariant algebra generators must have real coefficients. We may take the quotient over $a_v-\sum_{i} a_{i,k}=0$, then the quadratic function above corresponds to the Tits form of signature $(1,N(M-1))$. Furthermore, it is generally known that a hyperbolic Coxeter group $W_M$ is Zariski dense in the Lorentz group $O(1,N(M-1))$ of this form (see \cite{YH}). Therefore, all the invariants of the quotient symmetric algebra are generated by the quadratic function above. Now the original statement follows by taking the limit $M\rightarrow \infty$. $\square$\\

\subsection{Parameterization of reduced words.} Now we would like to consider the subdiagram $X$ of the star-shaped diagram with the central vertex removed. It is clear that the Weyl group $W_X$ of this subdiagram is isomorphic to $\prod_{i=1}^N S_{\infty,k}$ - the product of $N$ infinite symmetric groups. We wish to somehow parameterize the set $W^X$ of $X$-reduced elements of the Weyl group $W$ (see Bourbaki, Ch. 4, Ex. 3 \cite{Bou}).\\

\textbf{Proposition 4.2.} The $W$-orbit of the vector $\epsilon_v$ with $a_v=1,a_{i,k}=0$ consists of some vectors with integer coefficients $a_{i,k},a_v$, such that
\begin{equation}
    \forall k, \quad a_v-\sum_{i}a_{i,k}=1, \quad -(N-2)a_v^2+\sum_{i,k}a_{i,k}^2=-(N-2), \quad a_{i,k}\ge 0, \quad (N-2)|a_{i,k}.
\end{equation}
Moreover, the stabilizer of $\epsilon_v$ is precisely the subgroup $W_X$.\\

\textbf{Proof.} As we have already seen, the functions above are $W$-invariants, so they have to be constant on the orbit. For the divisibility property, note that $v$ has $(N-2)|a_{i,k}$ and the only generator that might change is $s_v$, but from the formula for the generator it is evident that if a vector $q$ satisfies this property then $s_v(q)$ does as well. As for the inequality, note that
\begin{equation}
    \epsilon_v=\alpha_v+\sum_{i\le M,k} \alpha_{i,k} -\sum_{k}\epsilon_{M+1,k}
\end{equation}
and $\alpha_v+\sum_{i\le M,k} \alpha_{i,k}$ is a positive root. If we choose any fixed element $w$ of the Weyl group $W$, then in its reduced form it must be finitely supported on the Dynkin diagram, therefore the set of inversion roots for $w$ must lie in a finitely generated root sub-lattice. On the other hand, we may choose $M$ arbitrarily large, so that $w$ does not affect $\epsilon_{M+1,k}$ and $\alpha_{M,k}$, therefore
\begin{equation}
    w(\alpha_v+\sum_{i\le M,k} \alpha_{i,k})
\end{equation}
must be a positive root with $a_{M+1,k}=1$. Additionally, the support of any root must be connected (Bourbaki, p. 173, Cor. 3 \cite{Bou}), so $a_v \neq 0$. Therefore from the description of the positive part of our Lie algebra from the previous chapter it is evident that $w(\epsilon_v)$ has nonnegative coefficients.\\

Finally, it is clear that the stabilizer of $\epsilon_v$ contains $W_{X}$. Suppose that $w\in W$ stabilizes $\epsilon_v$. Then in a matrix form (w.r.t. $a_v$ and $a_{i,k}$) the element $w$ looks like
\begin{equation}
    w=\begin{pmatrix}
        1 & v^t\\
        0 & A
    \end{pmatrix}.
\end{equation}

Note that 
\begin{equation}
    \begin{pmatrix}
        -(N-2) & -(N-2)v^t\\
        -(N-2)v & *
    \end{pmatrix}=
    w^t\cdot
    \begin{pmatrix}
        -(N-2) & 0\\
        0 & Id_{\infty}
    \end{pmatrix}
    \cdot
    w=
    \begin{pmatrix}
        -(N-2) & 0\\
        0 & Id_{\infty}
    \end{pmatrix},
\end{equation}
so $v=0$. Then $w$ must be an integer orthogonal matrix in the subspace defined by the equation $a_v=0$, but the only vectors with integer coefficients of length one are $\pm \epsilon_{i,k}$, so $w(\epsilon_{l,m})=\pm \epsilon_{i,k}$. By considering the invariants, we derive that $m=k$ and the sign is $+$, so $w$ lies in $W_X$. $\square$\\

\textbf{Remark 6.} The dual statement for the non-negativity of $w(\epsilon_v)$ would mean that all the elements of $\epsilon_{i,k}$-orbit are such that $a_v$ for them is non-positive.\\

From Proposition 4.2 it follows that all $X$-reduced (left or right) elements of the Weyl group are parameterized by the elements of the orbit $W\epsilon_v$. We may construct the following graph of states. Let us take $\epsilon_v$ and assign level $1$ to this state. Then we may act on this state by the Weyl generators. If we obtain a new state, let us assign it level 2 and connect the new state with the original one by an edge, and so on. It is clear that such graph might only have edges between states of level $l$ and $l+1$ or an edge between two states of the same level $l$ for some $l$. We may consider the algorithm of descent which takes the biggest positive element in a leg and pulls it closer to the central entry (this may happen in a non-unique way). Then if all $a_{1,k}$'s are the biggest elements in each leg the algorithm applies $s_v$ to this state, and so on.\\
\begin{center}
\begin{tikzpicture}
\node at (0,0) {$4$};
\node at (0,0.5) {$2$};
\node at (0,1) {$1$};
\node at (-0.3,-0.5) {$2$};
\node at (-0.6,-1) {$1$};
\node at (0.3,-0.5) {$1$};
\node at (0.6,-1) {$2$};

\draw[thick,->] (0.5,0) -- (1.5,0);

\node at (0 +2,0) {$4$};
\node at (0 +2,0.5) {$2$};
\node at (0 +2,1) {$1$};
\node at (-0.3 +2,-0.5) {$2$};
\node at (-0.6 +2,-1) {$1$};
\node at (0.3 +2,-0.5) {$2$};
\node at (0.6 +2,-1) {$1$};

\draw[thick,->] (0.5+2,0) -- (1.5+2,0);

\node at (0 +4,0) {$2$};
\node at (0 +4,0.5) {$0$};
\node at (0 +4,1) {$1$};
\node at (-0.3 +4,-0.5) {$0$};
\node at (-0.6 +4,-1) {$1$};
\node at (0.3 +4,-0.5) {$0$};
\node at (0.6 +4,-1) {$1$};

\draw[thick,->] (0.5+4,0) -- (1.5+4,0);

\node at (0 +6,0) {$2$};
\node at (0 +6,0.5) {$0$};
\node at (0 +6,1) {$1$};
\node at (-0.3 +6,-0.5) {$0$};
\node at (-0.6 +6,-1) {$1$};
\node at (0.3 +6,-0.5) {$1$};

\draw[thick,->] (0.5+6,0) -- (1.5+6,0);

\node at (0 +8,0) {$2$};
\node at (0 +8,0.5) {$0$};
\node at (0 +8,1) {$1$};
\node at (-0.3 +8,-0.5) {$1$};
\node at (0.3 +8,-0.5) {$1$};

\draw[thick,->] (0.5+8,0) -- (1.5+8,0);

\node at (0 +10,0) {$2$};
\node at (0 +10,0.5) {$1$};
\node at (-0.3 +10,-0.5) {$1$};
\node at (0.3 +10,-0.5) {$1$};

\draw[thick,->] (0.5+10,0) -- (1.5+10,0);

\node at (0 +12,0) {$1$};

\end{tikzpicture}\\
    \textbf{Fig. 2.} Example of the algorithm of descent for $N=3$. The length of the corresponding word is equal to $6$.
\end{center}
\vspace{0.5 cm}

\textbf{Proposition 4.3.} The algorithm of descent applied to an element of orbit $W\epsilon_v$ always terminates at $\epsilon_v$.\\

\textbf{Proof.} Assume without loss of generality that $a_{1,k}$ is the largest positive number in each individual leg. Let $x_{n,k}=a_{n,k}/a_{1,k}, n>1$, then $0\le x_{n,k}\le 1$ and
\begin{equation}
    -(N-2)a_v^2 +\sum_k a_{1,k}^2 (1+\sum_{n>1}x_{n,k}^2)=-(N-2),
\end{equation}
\begin{equation}
    a_v = 1+a_{1,k}(1+\sum_{n>1}x_{n,k}).
\end{equation}
From this we get
\begin{equation}
    a_{1,k}=\frac{a_v-1}{1+\sum_{n>1}x_{n,k}},
\end{equation}
so
\begin{equation}
    (N-2)(1-a_v^2)+\sum_{k} (a_v-1)a_{1,k}\frac{1+\sum_{n>1}x_{n,k}^2}{1+\sum_{n>1}x_{n,k}}=0
\end{equation}
and
\begin{equation}
    \sum_{k} a_{1,k}\frac{1+\sum_{n>1}x_{n,k}^2}{1+\sum_{n>1}x_{n,k}}=(N-2)(a_v+1).
\end{equation}
Therefore
\begin{equation}
    [(N-1)a_v-\sum_{k}a_{1,k}]-a_v=(N-2)a_v-\sum_{k}a_{1,k}=\sum_{k} a_{1,k}\frac{\sum_{n>1}x_{n,k}(x_{n,k}-1)}{1+\sum_{n>1}x_{n,k}}-(N-2)<0.
\end{equation}
This means that we may reduce $a_v$ until one of the maximal $a_{1,k}$ becomes $0$, but in that case $a_v=1$ and all $a_{i,k}=0$. $\square$\\

\textbf{Remark 7.} The $W$-orbit of $\epsilon_v$ does not coincide with the set of the integer points satisfying the conditions above: for $N=3$ consider the vector $\phi$ with
\begin{equation}
    a_{v}=4,a_{1,1}=a_{2,1}=a_{3,1}=a_{1,2}=a_{2,2}=a_{3,2}=1, a_{1,3}=-1,a_{2,3}=a_{3,3}=2,
\end{equation}
then $s_v(\phi)$ has
\begin{equation}
    a_v=7,a_{1,3}=a_{2,3}=a_{3,3}=2, a_{1,1}=a_{1,2}=4, a_{2,1}=a_{3,1}=a_{2,2}=a_{3,2}=1.
\end{equation}
In order to obtain a full description of the orbit, one has to describe the dual orbit of an $\epsilon_{n,k}$, which is an equally difficult task.\\

From the uniqueness of the $X$-reduced element of a coset in $W/W_X$ it follows that if two states at levels $m$ and $n$ are connected by two paths of length $|n-m|$, then the elements of the Weyl group corresponding to the products of the corresponding $|n-m|$ generators must be equal. \\

We will need the following Lemma.\\

\textbf{Lemma 4.4.} Let $s_{i_1}s_{i_2}\dots s_{i_n}s_v$ be a reduced decomposition of an $X$-reduced element of $W$, then the inversion root
\begin{equation}
    s_v s_{i_n} \dots s_{i_2}(\alpha_{i_1})
\end{equation}
has positive $a_v$.\\

\textbf{Proof.} We wish to prove that
\begin{equation}
    (\epsilon_v, s_v s_{i_n} \dots s_{i_2}(\alpha_{i_1}))<0,
\end{equation}
but firstly, let us note that the root
\begin{equation}
    s_vs_{i_{n}}\dots s_{i_2}(\alpha_{i_1})
\end{equation}
is positive, therefore if it has some nonzero $a_v$, then it must be positive. Now,
\begin{equation}
    (\epsilon_v, s_v s_{i_n} \dots s_{i_2}(\alpha_{i_1}))=(s_{i_2}\dots s_{i_n}s_v\epsilon_v, \alpha_{i_1}).
\end{equation}
The element $s_{i_2}\dots s_{i_n}s_v\epsilon_v$ is a previous state of the state $s_{i_1}s_{i_2}\dots s_{i_n}s_v\epsilon_v$ and the condition that
\begin{equation}
    (s_{i_2}\dots s_{i_n}s_v\epsilon_v, \alpha_{i_1}) \neq 0
\end{equation}
means precisely that $s_{i_2}\dots s_{i_n}s_v\epsilon_v$ and $s_{i_1}s_{i_2}\dots s_{i_n}s_v\epsilon_v$ are different states. $\square$\\

\textbf{Corollary 4.5.} From the positivity in the Lemma 4.4 we can see that the algorithm of descent never increases the level of a state.\\

It is possible to prove an even stronger statement.\\

\textbf{Proposition 4.6.} The graph of states does not have any edges on the same level.\\

\textbf{Proof.} Suppose that it does, then we have
\begin{equation}
    s_k s_{i_1}\dots s_{i_n}(\epsilon_v) = s_{j_1}\dots s_{j_n}(\epsilon_v)
\end{equation}
for some $(\emptyset, X)$-reduced words $s_{i_1}\dots s_{i_n}$ and $s_{j_1}\dots s_{j_n}$. Clearly $s_k s_{i_1}\dots s_{i_n}$ is reduced, otherwise $l(s_k s_{i_1}\dots s_{i_n})<n$ and the state $s_{j_1}\dots s_{j_n}(\epsilon_v)$ is actually of some level $l<n+1$ which is absurd. Then
\begin{equation}
    s_k s_{i_1}\dots s_{i_n} =s_{j_1}\dots s_{j_n} \cdot w_x,\quad  w_x\in W_X
\end{equation}
and $n+1=l(s_{j_1}\dots s_{j_n} \cdot w_x)=n+l(w_x)$, so $w_x=s_m$ for some $m$. The statement of the Proposition is then equivalent to the following: this situation can not occur for $s_{i_1}\dots s_{i_n} \neq s_{j_1}\dots s_{j_n}$. Indeed, suppose that $s_{i_1}\dots s_{i_n}$ is not $(\{k\},\o)$-reduced, then
\begin{equation}
    s_k s_{i_1}\dots s_{i_n} = s_k s_k s_{i_1'}\dots s_{i_{n-1}'}=s_{i_1'}\dots s_{i_{n-1}'},
\end{equation}
so the word $s_k s_{i_1}\dots s_{i_n}$ is not reduced, which is impossible. But then we definitely know that $s_{j_1}\dots s_{j_n}$ is not $(\{k\},\o)$-reduced and
\begin{equation}
    s_{j_1}\dots s_{j_n} = s_k s_{j_1'} \dots s_{j_{n-1}'}, 
\end{equation}
where $s_{j_1'} \dots s_{j_{n-1}'}$ is $(\{k\},X)$-reduced, so 
\begin{equation}
    s_{i_1}s_{i_2}\dots s_{i_n}=s_{j_1'} \dots s_{j_{n-1}'} s_{m},
\end{equation}
both words $s_{i_2}\dots s_{i_n}$ and $s_{j_1'} \dots s_{j_{n-1}'}$ are $(\emptyset,X)$-reduced and $s_{i_1}s_{i_2}\dots s_{i_n}$ is reduced, so we are in the same situation and we can induct on $n$. $\square$\\

\textbf{Corollary 4.7.} The algorithm of descent is the fastest algorithm in terms of number of steps (it takes exactly $l$ steps to get from a state of level $l+1$ to $\epsilon_v$).\\


\section{Weyl-Kac Formula for KM$(D,\textbf{t})$}
This section will follow closely the work of A. Pakharev \cite{AP}. Let $D$ be a general graph as in section 2.\\

\textbf{Definition 5.1.} The parabolic category $\mathcal{O}$ is the category of $\mathbb{Z}^{\oplus 2|\textbf{V}(D)|}$-graded $\text{KM}(D,\textbf{t})$-modules in Ind$\Rep(\prod_{e\in \textbf{E}(D)} GL_{t_e})$ which are finitely generated as $\mathfrak{u}_-$-modules and such that the action of $\mathfrak{u}_+$ is locally nilpotent.\\

The $\mathbb{Z}^{\oplus 2|\textbf{V}(D)|}$-grading here comes from the action of $\bbone_{v}$ and $\bbone_v'$ for $v\in \textbf{V}(D)$. Let $P_{\infty}$ be an infinitely generated abelian group of integer weights for graded objects from $\Rep(\prod_{e\in \textbf{E}(D)} GL_{t_e})$ spanned by symbols $\epsilon_{i,e},i\in \mathbb{Z}\setminus \{0\},e\in \textbf{E}(D)$ and $\epsilon_v,\epsilon_v'$ corresponding to the action of $\bbone_v,\bbone_v'$. Let $P_{\infty}^+$ be the subset of dominant weights in $P_{\infty}$ corresponding to finite sums
\begin{equation}
    \phi=\sum_{e\in \textbf{E}(D),i\in \mathbb{Z}\setminus \{0\}} a_{i,e} \epsilon_{i,e} + \sum_{v\in \textbf{V}(D)}a_{v}\epsilon_v +\sum_{v\in \textbf{V}(D)}a_{v}'\epsilon_v',
\end{equation}
\begin{equation}
    a_{i+1,e}\ge a_i \text{ if } i>0,\quad a_{i-1,e}\ge a_i \text{ if } i<0, \quad a_v' +\sum_{e\in E(D),s/t(e)=v} a_{\pm 1, e} \ge (N-2)a_v, \quad v\in \textbf{V}(D)
\end{equation}
with non-positive $a_{i,e}$, where the sign is $+$ if $v$ is the target of $e$ and $-$ if it's the source of $v$. The symbols $\epsilon_{i,e}$ for $i>0$ correspond to the weights of subobjects of $V_{e}^{\otimes n}$ and the symbols with $i<0$ correspond to the weights of subobjects of $(V_{e}^*)^{\otimes n}$. Let $\mathcal{C}^a$ be the category of $\mathbb{Z}^{\oplus 2|\textbf{V}(D)|}$-graded simple objects of $\Rep(\prod_{e\in \textbf{E}(D)} GL_{t_e})$ which are parameterized by $\phi\in P_{\infty}^+ \times \mathbb{Z}^{\oplus 2|\textbf{V}(D)|}$.\\

Due to the choice of the presentation of the Lie algebra $\mathfrak{gl}_{\infty}$ the highest weight of $V_e$ is $-\epsilon_{-1,e}$ and the highest weight of $V_e^*$ is $-\epsilon_{1,e}$. Therefore, to any highest weight $\phi$ we may associate a simple module $L_{\phi}\in \mathcal{C}^a$ which corresponds to the tensor product over $e\in \textbf{E}(D)$ of objects
\begin{equation}
    L_{e,[(-a_{-1,e},-a_{-2,e},\dots),(\dots, -a_{2,e},-a_{1,e})]}
\end{equation}
with an unchanged grading (via $a_v,a_v',v\in \textbf{V}(D)$).\\

\textbf{Definition 5.2.} For any simple object $L_{\phi}\in \mathcal{C}^a$, define the parabolic Verma module as $M(\phi)=U(\text{KM}(D,\textbf{t}))\otimes_{U(\mathfrak{l} \oplus \mathfrak{u}_+)} L_{\phi}$ with zero action of $\mathfrak{u}_+$ on $L_{\phi}$.\\

By using the standard argument it is easy to see that the module $M(\phi)$ admits the unique simple graded quotient $L(\phi)$. Now, consider a Coxeter group $W_D=\prod_{v\in \textbf{V}(D)} W_v$ where $W_v$ is the Weyl group of the vertex $v$ from the previous section with nontrivial action on $\epsilon_{i,e}$'s for $i>0$ if $t(e)=v$ or for $i<0$ if $s(e)=v$. The extended 
(w.r.t. $\epsilon_v'$) action of $W_v$ is described as the usual permutation action on $\epsilon_{i,e}$ by $s_{i,e}$, $s_{1,e}(\epsilon_v)=\epsilon_v$ and $s_{1,e}(\epsilon_v')=\epsilon_v'$. Also, for $e\in \textbf{E}(D)$, such that $s(e)=v$ or $t(e)=v$ we have
\begin{equation}
    s_{\alpha_v} (\epsilon_{\pm 1, e})=-\epsilon_v - \sum_{l\in \textbf{E}(D), s/t(l)=v,l\neq e} \epsilon_{\pm 1, l}, 
\end{equation}
\begin{equation}
    s_{\alpha_v} (\epsilon_v)= (N-1)\epsilon_v -(N-2)\sum_{e\in \textbf{E}(D), s/t(e)=v} \epsilon_{\pm 1, e},
\end{equation}
\begin{equation}
    s_{\alpha_v} (\epsilon_v')=\epsilon_v'-\epsilon_v - \sum_{e\in \textbf{E}(D), s/t(e)=v} \epsilon_{\pm 1, e}.
\end{equation}
Again, the sign choice is made in such a way that $W_v$ acts in a "vicinity"\: of the vertex $v$ (e.g. the sign in $\epsilon_{\pm 1, e}$ is $+$ if $v=t(e)$ and $-$ if $v=s(e)$). The action of $W_v$ on other symbols is trivial. Let $W_v^v$ be the subset of $D\setminus \{v\}$-reduced words of the Weyl group $W_v$.\\

\textbf{Lemma 5.3.} Consider a root system $R$ associated to a Dynkin diagram $D$ and let $X$ be a subset of $D$ with a finite Weyl group $W_X$, then an inversion root $\alpha\in R$ of a $(\emptyset, X)$-reduced word $w\in W$ satisfies supp$(\alpha)\cap D\setminus X \neq \emptyset$.\\

\textbf{Proof.} Let us take the longest element $w_X$ in $W_X$, then we have
\begin{equation}
    u = w \cdot w_X,\quad l(u) = l(w)+l(w_X)
\end{equation}
Note that due to the decomposition of $u$ the set $R_X^+$ of all positive roots supported at $X$ lies in the set of inversion roots of $u$, therefore $w$ must send $w_X R_X^+=R_X^-$ to the set of negative roots $R^-$. In particular, we have $w(R_X^+)\subset R^+$ which is equivalent to the original statement. $\square$\\

\textbf{Remark 8.} In general, the group $W_X$ does not have to be finite for the statement to be true.\\

\textbf{Theorem 5.4. (Kac-Weyl formulas)} Let $\phi\in P_{\infty}^+$ be a dominant weight and let
\begin{equation}
    w = (w_v)_{v\in \textbf{V}(D)}\in \prod_{v\in \textbf{V}(D)} W_v^v,
\end{equation}
then in the $K$-group of the $\mathbb{Z}^{2|\textbf{V}(D)|}$-graded Deligne category Ind$\Rep(\prod_{e\in \textbf{E}(D)}GL_{t_e})$ we have the following identity
\begin{equation} \label{eq:Character}
    ch(L(\phi))= \sum_{w\in \prod_{v\in \textbf{V}(D)} W_v^v} (-1)^{\sum_{v\in \textbf{V}(D)}l(w_v)} ch(M((\prod_{v\in \textbf{V}(D)}w_v)\cdot \phi)).
\end{equation}
Additionally, the formula
\begin{equation} \label{eq:Denominator}
    \frac{1}{ch(M(0))}= \sum_{w\in \prod_{v\in \textbf{V}(D)} W_v^v} (-1)^{\sum_{v\in \textbf{V}(D)}l(w_v)} ch(L((\prod_{v\in \textbf{V}(D)}w_v)\cdot 0))
\end{equation}
also holds true.\\

\textbf{Proof.} The RHS of the formula makes sense, because if $\phi$ is a dominant weight, then $w\cdot \phi,w\in W^X$ for a subdiagram $X$ of a Dynkin diagram is an $X$-dominant weight. Following the proof of Theorem 3.13 from A. Pakharev's work \cite{AP}, we can say that each $\mathbb{Z}^{2|\textbf{V}(D)|}$-graded component on the RHS and the LHS involves only finitely many irreducible objects from $\Rep(\prod_{e\in \textbf{E}(D)}GL_{t_e})$. Note that a fixed graded part of $L(\phi)$ for $\text{KM}(D,\textbf{t})$ and the part of $L(\phi)$ with the same grading for an associated Lie algebra with integer parameters $t_e$ have the same decomposition into simple objects for large enough $t_e$. Indeed, since every submodule of an associated $L(\phi)$ in the fixed grading has a highest weight with respect to $\prod_{e\in \textbf{E}(D)}GL_{t_e}$, in order for us to construct a maximal quotient $L(\phi)$ of $M(\phi)$, it is necessary and sufficient to find submodules annihilated by the adjoint action of $\mathfrak{u}_+$. However, it is clear that the annihilation condition can be checked in some finite filtration degree of $\mathfrak{u}_+$. This argument shows that fixed degree parts of characters of $L(\phi)$ coincide for large enough $t_e$.\\

The formulas on the RHS for finite case and a categorical one involve different sets: $W^X$ where $X$ is the set of edge vertices and the product of $W_v^v$'s. However, due to Lemma 4.4 (infinite case for a vertex with $N\ge 3$), \cite{AP} (cases $N=1$ and $N=2$ are trivial) and Lemma 5.3 (note that $W_X$ in this case is $\prod_{e\in \textbf{E}(D)} S_{t_e}$ which is finite) we know that both formulas involve only finitely many terms from both sets for a fixed graded part of the RHS. So if we consider large enough $t_e$'s then all words with small length from $W^X$ decompose into products of small and independent commuting words from $W_v^v$'s. It is also evident that for large $t_e$ the RHS for both cases involve the same irreducible objects (for a fixed graded part as well). Therefore, the two equalities above and the equality in the Kac-Weyl formula imply the first statement. The second statement follows from the first just like in Corollary 3.14 in \cite{AP} (this is where the transcendence condition is used). $\square$\\

\textbf{Remark 9.} The Kac-Weyl character formula in this Theorem has an interesting property. Namely, the summation on the RHS is taken over non-interacting pieces of data $W_v^v$ near each individual vertex $v$.\\

This Theorem and the descent algorithm from the previous chapter allows us to compute a partial character of any simple integrable $\text{KM}(D,\textbf{t})$-module $L(\phi)$ within some finite region of the $\mathbb{Z}^{2|\textbf{V}(D)|}$-lattice. However, we could not find good formulas for description of $N$-tuples of partitions arising near any vertex of $D$ of valency $N\ge 3$. \\

\textbf{Example 5.5.} Consider the setting of section 4 with $N=3$. In this case we have only one vertex $v$ to work with. Let us write out the denominator formula \eqref{eq:Denominator}. For this example we write weights $\phi$ by using the following coordinate expression:
\begin{equation}
    \phi = (a_v; a_{-1,1},a_{-2,1},\dots ;a_{-1,2},a_{-2,2},\dots; a_{-1,3},a_{-2,3},\dots).
\end{equation}
Note that here we don't have to deal with coefficients $a_{i,k},i>0,1\le k \le 3$, because they will always be zero. Similarly we can omit the coordinate $a_v'$. For any reduced word $w\in W_v^v$ we have $w\cdot 0 = -\sum_{\alpha\in I(w)} \alpha$ where $I(w)$ is the set of inversion roots for $w$. Let us compute the denominator formula up to the level $4$. At this level we have $6$ reduced words
\begin{equation}
    s_{2,k}s_{1,k}s_v, \quad s_{1,l}s_{1,p}s_v, \quad 1\le k\le 3, \quad 1\le l<p \le 3.
\end{equation}
We also have $3$ words at level $3$
\begin{equation}
    s_{1,k}s_v, \quad 1\le k\le 3
\end{equation}
and one word $e$ and $s_v$ at levels $1$ and $2$ correspondingly. Then the denominator formula reads
\begin{multline}
    \frac{1}{ch(M(0))} = ch(L(0))- ch(L(-1;-1,0,\dots;-1,0,\dots;-1,0,\dots))+\\+ ch(L(-2;-2,0, \dots; -2,0,\dots; -1,-1,0,\dots))+ch(L(-2;-2,0, \dots; -1,-1,0,\dots; -2,0,\dots))+\\+ch(L(-2; -1,-1,0,\dots; -2,0, \dots; -2,0,\dots))-ch(L(-3;-3,0, \dots; -2,-1,0,\dots; -2,-1,0,\dots))-\\-ch(L(-3; -2,-1,0,\dots;-3,0, \dots; -2,-1,0,\dots))-\\-ch(L(-3; -2,-1,0,\dots; -2,-1,0,\dots;-3,0, \dots))-\\-ch(L(-3;-3,0, \dots; -3,0,\dots; -1,-1,-1, 0,\dots))-\\-ch(L(-3;-3,0, \dots; -1,-1,-1, 0,\dots; -3,0,\dots))-\\-ch(L(-3; -1,-1,-1, 0,\dots; -3,0, \dots;  -3,0,\dots)) +\dots
\end{multline}
    
\selectlanguage{english}

\end{document}